\input amstex
\magnification=1200
\documentstyle{amsppt}
\topmatter
\rightheadtext{~}
\leftheadtext{~}

Fields Inst.\ Communications,
{\bf 24} (1999) 311--323

ITEP-TH-63/98

\vskip -.5cm

\title
Symplectic geometry on moduli spaces of holomorphic bundles over complex
surfaces
\endtitle

\author Boris Khesin and Alexei Rosly \endauthor

\date{October 1998}\enddate
\vskip 1in

\abstract

We give a
comparative description of the Poisson structures on the moduli
spaces of flat connections on {\it real}  surfaces and holomorphic
Poisson structures on the moduli spaces of holomorphic bundles on
{\it complex} surfaces.  The symplectic leaves of the latter are
classified by restrictions of the bundles to certain divisors.
This can be regarded as fixing a ``complex analogue of
the holonomy" of a connection along a ``complex analogue of the
boundary" in analogy with the real case.
\endabstract

\endtopmatter
\document
\widestnumber\key{EFWS}
\baselineskip 14pt
\NoBlackBoxes

\head {  Introduction}
\endhead

In this  note we discuss the geometry of the momentum map  for
gauge groups
in the following two cases with the aim of emphasizing
the analogy between them. We start by recalling a
description of the Poisson structure on the moduli space of flat
connections
over a two-dimensional real surface. Our main interest  is related to
a description of analogous
structures on moduli of holomorphic bundles over a two-dimensional
complex
surface. In the first case we deal with smooth objects while in the
second one -- with complex analytic objects.
Our interest in this subject comes mainly from a desire to understand
the origin of a symplectic structure (refs. \cite{Mu, Ko}) and of
a Poisson structure (refs. \cite{Bon, Bot})
on moduli of holomorphic bundles over complex surfaces in a way
which would be parallel to the consideration of flat
connections over real surfaces. It is worth mentioning that the both
cases above are also of interest from the mathematical physics point of
view.  Other related important results on the
geometry of moduli spaces of holomorphic bundles on certain complex
surfaces deal with the study of
the symplectic structure related to a K\"ahler form, or
a hyperk\"ahler structure (cf., e.g., refs. \cite{Do, KN, LMNS}).

We adopt the viewpoint  of considering
 the case of holomorphic
bundles on complex surfaces as a certain
complexification of the case of flat connections on
real surfaces. This approach is parallel to the geometric complexification
method  suggested by V.~Arnold in ref. \cite{Ar}.
To be more precise, rather than formally complexifying,
we replace locally constant sheaves (corresponding to flat connections)
by sheaves of holomorphic sections
(or, one could say, $d/dx$ is replaced by $\partial/\partial\bar z$).
Being rather simple by itself this leads however to
certain curious ideas some of which we will sketch below
(cf., also refs.\cite{Wi, EF, FK, Kh, DT, Th, FKT, KR}).

The consideration of flat connections rests on the notions of holonomy
and curvature. One needs to equate the curvature  to zero, which leads
to the ``flatness" condition, while
the holonomy is crucial as the ``only remaining" part of
what can characterize a
flat connection  modulo gauge  transformations.
Thus, to pass to ``complexified" objects, i.e., to holomorphic bundles over
complex surfaces, we need to know what  the complex analogues of
holonomy and curvature are. It is somewhat easier with curvature.
For a (0,1)-connection (i.e.,
$\bar\partial$-connection) one can define its curvature (0,2)-form.
However, as we shall see below (cf., \S 3), a better analogue
of the curvature (in the context of
symplectic geometry on the space of $\bar \partial$-connections) will be
a  certain (2,2)-form. This form is the
wedge product of the above (0,2)-form with a  meromorphic (or holomorphic)
(2,0)-form on  the surface.  The meromorphic
(2,0)-form will play the
role of orientation of the surface in the complex analytic situation.
We specifically concentrate on the case of
a meromorphic (2,0)-form with  logarithmic singularities. This gives
rise to a complex analogue of the notion of a real surface with
boundary, namely,   a complex surface with the
polar set of a (2,0)-form.

The problem of finding a proper complex analogue of the notion of
holonomy of a flat connection (or, better to say, monodromy)
is more intricate. Let us
first consider the case of a real surface with boundary. Then, for a
flat connection on it, we have certain group elements, the monodromies,
associated with every closed
loop in the surface. In the case when the loop is
homotopic to a boundary component we say that we deal with a monodromy
around a hole in the surface. The monodromies around the holes play a
distinctive role in the study of Poisson geometry of the moduli space of
flat connections. The latter space is, in fact, a Poisson manifold whose
Poisson structure is degenerate in the case when the surface has
holes (i.e., non-empty boundary). The symplectic leaves of that Poisson
manifold can be defined by fixing the conjugacy classes of the
monodromies
around the holes. Thus one can find out a proper ``complexified" notion
of the monodromy around a hole provided one is able to single out the
symplectic leaves in the Poisson manifold of moduli of holomorphic
bundles on
a complex surface with ``boundary" (in the sense mentioned above). This
will be discussed  in \S 3.
In the case of a loop which is not homotopic to the
boundary we do not know yet a proper complex analogue of the monodromy.
It would be interesting to see what should replace ``a loop," a concept
from homotopy theory, in the complex analytic setting.

It turns out, however, that at least the corresponding homology theory
can be constructed. The above approach leads one to
certain complex analytic analogues of the notions of
chains, boundary, and cycles. This is mentioned in \S 2.

\bigskip

\subhead \S 1. Real case: Poisson structures on moduli spaces of flat
connections
\endsubhead

\medskip

First we recall several results on Poisson structures
on the space of flat connections on {\it real} surfaces with boundary.
Our consideration of holomorphic bundles on {\it complex} surfaces
below,
in \S 3, will be parallel to the case of real surfaces.
We follow the papers [AB, FR] in the exposition of the real case.

In  the  real case $G$ stands for a simple simply connected {\it
compact} Lie group, and  $\goth g=\operatorname{Lie} G$ is its Lie
algebra
with a chosen nondegenerate invariant quadratic form, which we denote by
$\operatorname{tr}$.
Let $\Sigma$  be an oriented compact surface which may have boundary
$\Gamma=\partial \Sigma$
consisting of several components, $\Gamma=\cup^k_1 \Gamma_j$.
Denote by $E$ a (trivial) principle $G$-bundle over $\Sigma$.
Let $\Cal A^\Sigma$ be the affine space of all smooth
connections in $E$.
It is convenient to fix any trivialization of $E$ and identify
$\Cal A^\Sigma$ with the vector space $\Omega^1(\Sigma,\goth g)$ of
smooth $\goth g$-valued 1-forms on the surface:
$$
\Cal A^\Sigma=\{ d+A~|~A\in\Omega^1(\Sigma,\goth g)\}~.
$$

\definition{Definition 1.1}
The space $\Cal A^\Sigma$ is in a natural way a symplectic manifold
with the symplectic structure
$$
W:=\int_\Sigma\operatorname{tr}(\delta A\wedge \delta A)~,
\tag{1}
$$
where $\delta$ is the exterior differential on $\Cal A^\Sigma$, and
$\wedge$ stands to denote the wedge product both on $\Cal A^\Sigma$
and $\Sigma$.
\enddefinition

\proclaim{Proposition 1.2}
The symplectic structure $W$ is invariant with
respect to the gauge transformations
$$
A\mapsto g^{-1}Ag+g^{-1}dg~,
$$
where $g$ is an element of the group of gauge transformations,
$G^\Sigma$, i.e., it is a smooth $G$-valued function on the surface
$\Sigma$.
\endproclaim

The infinitesimal gauge transformations forming a Lie algebra
$\goth g^\Sigma$ are generated on the symplectic manifold $\Cal
A^\Sigma$ by
certain Hamiltonian functions.

\proclaim{Proposition 1.3}
An infinitesimal gauge transformation $\epsilon$ is generated
by the Hamiltonian function
$$
H_\epsilon=\int_\Sigma\operatorname{tr}(\epsilon(dA+A\wedge A))-
\int_{\partial\Sigma}\operatorname{tr}(\epsilon A)~.
$$
\endproclaim

\demo{Proof}
Hamiltonian vector field $X$ corresponding to any Hamiltonian $H$
is defined by its action on functions $f(A)$ :
$$
L_Xf=\{H,f\}=\int_\Sigma\operatorname{tr}
\left(\frac{\delta H}{\delta A}\wedge\frac{\delta f}{\delta A}\right),
$$
where the latter expression is the Poisson bracket corresponding to eq.(1).
It suffices to consider the coordinate function $f(A)=A$ :
$$
L_XA=\{H,A\}=\frac{\delta H}{\delta A}.
$$
Then for the above Hamiltonian $H_\epsilon$ we have
$$
L_{X_\epsilon}A=\frac{\delta H_\epsilon}{\delta A}=\nabla\!_A\,\epsilon,
$$
where $\nabla\!_A\,\epsilon=d\epsilon+[A,\epsilon]$ is an infinitesimal
gauge transformation. Indeed, for $F(A):=dA+A\wedge A$ we have
$\delta F=\nabla\!_A\,\delta  A$, \,\, and then
$$
\align
\delta H_\epsilon=&\int_\Sigma\operatorname{tr}
\left(\epsilon\,\delta F\right)
 -\int_{\partial\Sigma}\operatorname{tr}(\epsilon\,\delta  A)\\
=&\int_\Sigma\operatorname{tr}
\left(\epsilon\,\nabla\!_A\,\delta  A\right)
 -\int_{\partial\Sigma}\operatorname{tr}(\epsilon\,\delta  A)
=\int_\Sigma\operatorname{tr}
\left(\delta  A\wedge\nabla\!_A\,\epsilon\right).
\endalign
$$
In the last equality we used the Stokes formula.\qed
\enddemo

The Hamiltonian function generating a given gauge transformation
is defined only up to an additive constant. Hence, generally
speaking, the Poisson
bracket between two such Hamiltonians reproduces the commutation
relation in the gauge algebra $\goth g^\Sigma$ only up to a cocycle:
$$
\{H_{\epsilon_1},H_{\epsilon_2}\}= H_{[\epsilon_1, \epsilon_2]}
+c(\epsilon_1, \epsilon_2)~.
$$

\proclaim{Proposition 1.4}
For the above choice of Hamiltonians
the cocycle is
$$
c(\epsilon_1, \epsilon_2)=\int_{\partial\Sigma}
\operatorname{tr}(\epsilon_1d\epsilon_2)~.\tag{2}
$$
\endproclaim

One can show that this well-known cocycle is nontrivial. Therefore one
can
define the momentum mapping not for the algebra of gauge
transformations, but
only for its central extension by the 2-cocycle ({2}).

To define this mapping we need some notations.
Let $\hat{\goth g}^\Sigma$ denote the Lie algebra of gauge
transformations centrally extended by the cocycle ({2}), and
$\hat{G}^\Sigma$ be the corresponding group.
The infinite-dimensional space $\hat{\goth g}^\Sigma$ is the space of
pairs
$(\epsilon, z)$, where $\epsilon$ is a $\goth g$-valued function on the
surface $\Sigma$ and $z$ is a real number.

We define  the space $(\hat{\goth g}^\Sigma)^*$,
dual to $\hat{\goth g}^\Sigma$, as consisting of triples $(F,C,x)$,
where $F$ is a $\goth g$-valued
2-form on $\Sigma$, $C$ is a $\goth g$-valued
1-form on the boundary of $\Sigma$, and $x$ is a real number.
The nondegenerate pairing $<\,,>$ between the spaces
$\hat{\goth g}^\Sigma$ and
$(\hat{\goth g}^\Sigma)^*$ is the following:
$$
<(F,C,x), (\epsilon, z)>=\int_\Sigma \operatorname{tr}(\epsilon F)-
\int_{\partial\Sigma} \operatorname{tr}(\epsilon C) +zx~.
$$

Let us consider the action of $\hat{G}^\Sigma$ on $\Cal A^\Sigma$
generated
by the action of $G^\Sigma$. That is to say, the center of
$\hat{G}^\Sigma$
acts trivially.

\proclaim{Proposition 1.5}
The centrally  extended group  $\hat G^\Sigma$ of gauge
transformations acts on $\Cal A^\Sigma$  in a Hamiltonian  way.  The
momentum map for the action of the corresponding gauge algebra
$\hat{\goth
g}^\Sigma$ is the mapping $\Cal A^\Sigma\to (\hat{\goth g}^\Sigma)^*$
given
by the curvature and by the restriction of the connection form to the
boundary:  $$ A\mapsto (dA+A\wedge A, ~~A|_{\partial\Sigma}, ~~1)~.  $$
\endproclaim

Let us introduce the following notation.
For a manifold $\Sigma$ and its submanifold $\Gamma\subset \Sigma$
denote by $G^\Sigma_\Gamma$ the group of gauge transformations on
$\Sigma$
``based on $\Gamma$":
$G^\Sigma_\Gamma=\{g\in C^\infty(\Sigma,G)~|~g|_\Gamma
=\operatorname{id}\}$, and by ${\goth g}^\Sigma_\Gamma$
the corresponding Lie algebra.

Modifying slightly the last proposition one gets the following

\proclaim{Corollary 1.6}
The group $G^\Sigma_\Gamma$ acts on $\Cal A^\Sigma$ in a Hamiltonian
way.
The momentum map for the action of the corresponding Lie algebra
${\goth g}^\Sigma_\Gamma$ is the mapping
$\Cal A^\Sigma\to (\goth g^\Sigma_\Gamma)^*$ given by the curvature:
$$
A\mapsto dA+A\wedge A~.
$$
\endproclaim

\remark{Remark} Note that the group $G^\Sigma_\Gamma$ is not centrally
extended, but still $G^\Sigma_\Gamma\subset \hat{G}^\Sigma$.
\endremark

Now consider the Hamiltonian reduction $\Cal
A^\Sigma//G^\Sigma_\Gamma$ of the space of connections $\Cal A^\Sigma$
with
respect to the group $G^\Sigma_\Gamma$ of gauge transformations equal to
the identity on the boundary $\Gamma=\partial \Sigma$.
This yields the space of flat
connections on $\Sigma$ modulo  gauge transformations from
$G^\Sigma_\Gamma$,
$$
\Cal M_{\Sigma,\Gamma}=\{d+A\in \Cal A^\Sigma\,|\,dA+A\wedge
A=0\}/G^\Sigma_\Gamma~.
$$

By definition of Hamiltonian reduction, the space $\Cal
M_{\Sigma,\Gamma}$
is symplectic (though, certainly, infinite-dimensional and,
generally speaking, with singularities).
It can be mapped to certain familiar Poisson manifolds.
It is well known that the space of
$G$-connections on a circle can be identified with the space
of coadjoint representation of the affine Kac--Moody algebra equipped
with the standard Kirillov--Kostant Poisson structure. The relation with
the Poisson (in fact, symplectic) structure on $\Cal M_{\Sigma,\Gamma}$ is
given by the following proposition.

\proclaim{Proposition 1.7}
The mapping from the space
$\Cal M_{\Sigma,\Gamma}$ to the Kac--Moody  coadjoint representation
space
sending a flat connection on the surface $\Sigma$ to its
restriction to a boundary component is a Poisson mapping.
\endproclaim

\demo{Proof} This mapping is essentially the momentum mapping for the
action of gauge transformations on the boundary. \qed
\enddemo

\remark{Remarks} {\it i)} Here and below we always mean the
 nonsingular parts of the moduli spaces when describing the
symplectic (or Poisson) structures on them.

{\it ii)} We are going to describe now the quotient of our symplectic
manifold by the Hamiltonian action of a group. The result is always a
{\it
Poisson} manifold, while the momentum map helps to determine the
symplectic
leaves in it (see, ref. [We]).
\endremark

Consider the quotient of the space $\Cal M_{\Sigma,\Gamma}$ by the
whole group $\hat{G}^\Sigma$ of centrally extended gauge
transformations.
The latter group acts on $\Cal M_{\Sigma,\Gamma}$ since gauge
transformations equal to the identity on the boundary form a
normal subgroup $G^\Sigma_\Gamma$ in $ \hat{G}^\Sigma$ .
The quotient space
$$
\Cal M_\Sigma=\{d+A\in \Cal A^\Sigma\,|\,dA+A\wedge
A=0\}/\hat{G}^\Sigma
$$
is a finite-dimensional Poisson manifold (with singularities).

\proclaim{Proposition 1.8} The space $\Cal M_\Sigma$
of flat $G$-connections modulo  gauge transformations
on a surface $\Sigma$ with holes inherits a Poisson structure from
the space of all (smooth) $G$-connections.
The symplectic leaves of this structure are parameterized by the
conjugacy classes of holonomies around the holes (that is, a symplectic
leaf is singled out by fixing the conjugacy class of the
holonomy around each hole).
\endproclaim

\demo{Proof} The symplectic leaves of $\Cal M_\Sigma$ are in one-to-one
correspondence with the coadjoint orbits of the (centrally
extended) affine Lie algebra on a circle (or the
direct sum of several copies of the
affine algebras if the boundary of the surface $\Sigma$ consists of
several components). These coadjoint orbits,  in turn,
are parameterized by the conjugacy classes of holonomies
around the circle. \qed
\enddemo

\remark{Remark}
The above proposition should not be understood as that the conjugacy
classes
of holonomies around the holes can be taken arbitrary, since the
holonomies of
a flat connection on the surface obey certain relations coming from the
fundamental group $\pi_1(\Sigma)$. For example, if $\Sigma$ is a sphere
with
$n$ holes then the product of all $n$ holonomies has to be
$\operatorname{id}\in G$  (provided
one has chosen the same base point and a convenient orientation for all
$n$
loops encircling the holes).
\endremark

\bigskip

\bigskip

\subhead \S 2. The Stokes--Leray formula \endsubhead
\medskip

In order to develop the symplectic geometry related with holomorphic
bundles on complex surfaces in a way analogous to what have been
considered in the last section we will need a complex analogue
of the Stokes formula. This
will be nothing but a simple multidimensional generalization of
the Cauchy formula.

{\bf  Higher-dimensional residue.}
Let $\gamma$ be a meromorphic $n$-form on
a  compact complex $n$-dimensional
manifold $M$ with poles  on a smooth complex
hypersurface  $N\subset M$. Here and below we consider the forms
with logarithmic singularities only (i.e., with the  first order
poles).\footnote{In this paper we restrict
ourselves to the case of top-degree meromorphic forms having singularities
on smooth divisors. In this case a logarithmic singularity is the same
as a first order pole.  It is the formulation
``$\gamma$ with logarithmic singularities'' which should be kept if
one would like to consider the situation of a non-smooth divisor of
poles.}
  Let $f$ be a function defining $N$ in a neighborhood of some
point $p\in N$.  Then locally, in a certain neighborhood $U(p)$, the
$n$-form
$\gamma$ can be decomposed into the sum $$ \gamma = \frac{df}{f}
\wedge\alpha
+ \beta ~, $$ where $\alpha$ and $\beta$ are holomorphic in $U(p)$. One
can
show, that the restriction $\alpha |_N$ is a well-defined (i.e.
independent
of $f$) holomorphic $(n-1)$-form on $N$.

\definition{Definition 2.1}
The holomorphic $(n-1)$-form $\alpha |_N$ on $N$ is called the residue
of the meromorphic form $\gamma$ and is denoted by
$~\text{res}_N\,\gamma$.
\enddefinition

\proclaim{Proposition 2.2}
Let $M,N$ and $\gamma$ be as above, and let $u$ be a smooth $(n-1)$-form
on $M$. Then the form $\gamma\wedge du$ is integrable on $M$, and
$$
\int_M \gamma\wedge du = 2\pi i\int_N \operatorname{res}_N\,
\gamma\wedge u~.
\tag{3}
$$
\endproclaim

\remark{Remarks}
{\it i)} The formula ({3}) is proved by applying the Stokes
formula to reduce the integral to the tubular neighborhood of $N$, and
then
by using the standard
Cauchy formula in the transversal direction to $N$ (see, e.g., ref. \cite{GS}).

{\it ii)} This relation can be, of course, generalized to the case
when $N$ is a normal crossing divisor in $M$ by modifying the above
definition
of a residue. Then, in particular, $\text{res}~\gamma$ will define
{\it meromorphic} forms on smooth components of $N$, the residues of
which (i.e. residues of residues) will sum to zero at the intersections
of components.

{\it iii)} We call this formula the Stokes--Leray formula for it
is a part of a much broader, than explained here, Leray theory (cf.,
ref. \cite{Le}), while, on the other hand, we are going to exploit it as a
complex analogue of the usual Stokes formula.

{\it iv)} Of course, it is only the $(0,n-1)$-part of the form
$u$ which is essential in eq.({3}). In the same way, we can
write $\bar\partial u$ instead of $du$:
$$ \int_M \gamma\wedge\bar\partial u =
  2\pi i\int_N \text{res}_N\,\gamma\wedge u~.
\tag{$3'$}
$$
\endremark

{\bf Digression on homology and cohomology.}
The last remark leads one to a simple idea of considering a pairing
of Dolbeault cochains with  certain ``geometric chains",
and thus to  the construction  of the corresponding homology
theory, which we will describe here roughly (and in a full detail in
ref. \cite{KR}).

Let $v$ be a smooth $(0,k)$-form on a complex manifold $M$
and let $(X,\alpha)$ be a pair consisting of a smooth $k$-dimensional
complex
submanifold $X$ in $M$ and a meromorphic $k$-form $\alpha$ on $X$
possessing logarithmic singularities on a smooth hypersurface $Y$ in
$X$. Now
let us consider the pairing between $(0,k)$-forms and the set of such
pairs given by the integral
$$
\int_X ~\alpha\wedge v~.
\tag{4}
$$
Note that the
meromorphic top-degree form $\alpha$ on $X$ is
the data which  allow us to integrate
$(0,k)$-forms  over submanifolds of {\it complex}
dimension $k$,  (i.e., to integrate $v$ over $X$). Therefore, this
meromorphic form $\alpha$  can be regarded as a holomorphic analogue of
{\it orientation} of the submanifold $X$.
Furthermore, if $v=\bar\partial u$,
then, by use of eq.($3'$),  the
integral over $X$ is reduced to the integral over its submanifold $Y$ of
one complex dimension less:
$$
\int_X \alpha\wedge\bar\partial u
= 2\pi i\int_Y  \text{res}_Y\,\alpha\wedge u~.  \tag{$3''$}
$$
 Thus we can speak of the pair
$(Y,~\text{res}_Y\,\alpha)$ as  a holomorphic analogue of the boun\-da\-ry
for the pair $(X,\alpha)$. In this sense, eq.($3''$) can be
viewed  on  as a complex analogue of the Stokes formula.
Suppose now that
the pair $(X,\alpha)$ is an analogue of a closed manifold, i.e.,
that $\alpha$ is holomorphic.
Then the integral ({4}) for a $\bar\partial$-closed form
$v$
depends in fact only on the Dolbeault cohomology class of $v$. This line
of reasoning can be developed to a homology theory (ref. \cite{KR}) which
plays the same role with respect to Dolbeault cohomology as
singular homology plays with respect to De Rham cohomology.

\remark{Remark} Thus, a pair $(M,\gamma)$ with a {\it holomorphic}
form $\gamma$ of
degree equal to $\text{dim}_{\Bbb C}\,M$ can be regarded as a
holomorphic
analogue of an {\it oriented closed} manifold. Sometimes it is
necessary to require that $\gamma$ {\it has no zeros} (then $M$ has to
be a
Calabi--Yau or an abelian manifold) in which case one can speak of a
holomorphic analogue of a {\it smooth oriented closed} manifold.
If $\gamma$
is meromorphic (rather than holomorphic) with only first order poles one
can speak of ``a manifold with boundary" (and the above remark about the
possible
zeros of $\gamma$ applies in this case as well).
\endremark

In the next section  we would like to exploit
the above understanding of what are the proper holomorphic analogues
of orientation and boundary
(a similar point of view was useful in
refs. \cite{FK, DT} for other gauge-theoretic constructions).

\medskip

\bigskip

\subhead \S 3. Complex case: Poisson structures on moduli spaces
of holomorphic bundles
\endsubhead
\medskip

Let $S$ be a compact {\it complex} surface ($\text{dim}_{\Bbb C}\,S=2$).
We are
going to describe a Poisson structure on the moduli space of holomorphic
vector bundles \footnote{~By the moduli space we shall always understand
a local universal family near a smooth point.} with a {\it complex}
reductive
group $G$ as the structure group ($G\subset GL(n,\Bbb C))$ on $S$. Let
us do
this in analogy with the consideration of flat connections in \S 1.
First of
all, in order to define an analogue of the symplectic structure in
eq.(1) we have to fix a holomorphic analogue of the orientation.
According to the heuristic argument in \S 2, we have to choose a
meromorphic
2-form on the surface $S$.  Let $\sigma$ be a meromorphic 2-form on $S$,
such that its divisor of poles $P$ is a smooth curve in $S$ and that
$\sigma$ has there  a logarithmic singularity.
 The curve $P$ will play the
role of the boundary of the surface in our considerations.  Let us
assume
additionally that $\sigma$ has no zeros (the situation analogous to a
smooth oriented
real surface).  Then $P$ is an anticanonical divisor on $S$ and
has to be an elliptic curve, or, may be, a number of nonintersecting
elliptic
curves.  These are analogous to the circles constituting the boundary of
a real smooth surface. In what follows we shall assume that $S$ is endowed
with such a 2-form $\sigma$. (Example: $S=\Bbb CP^2$ with a smooth cubic
as an anticanonical divisor. As a matter of  fact, many Fano surfaces fall
into this class.)
If it happens that $\sigma$ has no zeros and no poles (i.e., $S$
is ``oriented, without boundary") it means that we deal with either a K3
or an abelian surface.  (Note that the further consideration can be
extended with minimal changes to the case of a non-smooth divisor $P$,
in particular, to $P$ consisting of several
components intersecting transversally.
Example: $S=\Bbb CP^2$ with $\sigma=dxdy/xy$.)

Let $E$ be a smooth vector $G$-bundle over $S$ which can be endowed with
a holomorphic structure and $\text{End}\,E$ be the corresponding bundle of
endomorphisms with the fiber $\goth g=\text{Lie}(G)$. Let $\Cal A^S$
denote the infinite-dimensional affine
space of smooth $\bar\partial$-connections in $E$. By choosing a
reference holomorphic structure $\bar\partial_{\,0}$,
$\bar\partial_{\,0}^{\;2}=0$, in
$E$, the space $\Cal A^S$ can be identified with the vector space
$\Omega^{(0,1)}(S,\text{End}\,E)$ of $(\text{End}\,E)$-valued
$(0,1)$-forms on
$S$, i.e.
$$
\Cal A^S = \{\bar\partial_{\,0} +A\,|\, A\in
\Omega^{(0,1)}(S,\text{End}\,E)\}~.
$$
In what follows,
instead of $\bar\partial_{\,0}$, we shall write simply
$\bar\partial$ keeping in mind that this corresponds to a reference
holomorphic
structure in $E$ when it applies to sections of $E$ or associated
bundles.

\definition{Definition 3.1}
The space $\Cal A^S$ possesses a natural holomorphic symplectic
structure
$$
 W_{\Bbb C} := \int_S \sigma \wedge\operatorname{tr}(\delta A_1\wedge
\delta
A_2)~,
$$
where $\sigma$ is the holomorphic ``orientation" of $S$, while the other
notations are essentially the same as in Definition 1.1 above.
\enddefinition

After such a definition we can repeat the contents of  \S 1 more
or less word by word.

\proclaim{Proposition 3.2}
The symplectic structure $W_{\Bbb C}$ is invariant with
respect to the gauge transformations
$$
A\mapsto g^{-1}Ag+g^{-1}\bar\partial g~,
$$
where $g$ is an element of the group of gauge transformations,
i.e., the group of automorphisms of the smooth bundle $E$. Abusing
notation we denote this group by $G^S$.
\endproclaim

The infinitesimal gauge transformations forming the Lie algebra
$\goth g^S=\Gamma(S,\text{End}\,E)$ (where $\Gamma$ denotes the space of
$C^\infty$-sections) are generated on the symplectic manifold $\Cal A^S$
by certain Hamiltonian functions.

\proclaim{Proposition 3.3}
An infinitesimal gauge transformation $\epsilon$ is generated
by the Hamiltonian function
$$
H_\epsilon=\int_S \sigma \wedge \operatorname{tr}(\epsilon\,
(\bar\partial A+A\wedge A))-
2\pi i\int_P \operatorname{res}_P\,\sigma \wedge
\operatorname{tr}(\epsilon A)~.
$$
\endproclaim

\demo{Proof} It is the proof of this Proposition where the
Stokes--Leray
formula of  \S 2 is used instead of the usual Stokes formula being
the only modification in comparison with the proof of Proposition 1.3
in  \S 1. Now we have
$$
\align
\delta H_\epsilon=&\int_S\sigma \wedge \operatorname{tr}
\left(\epsilon\,\delta F\right)
 -\int_P \operatorname{res}_P\,\sigma \wedge
\operatorname{tr}(\epsilon\,\delta  A)\\
=&\int_S\sigma \wedge \operatorname{tr}
\left(\epsilon\,\nabla\!_A\,\delta  A\right)
 -\int_P \operatorname{res}_P\,\sigma \wedge
\operatorname{tr}(\epsilon\,\delta  A)
=\int_S\sigma \wedge \operatorname{tr}
\left(\delta  A\wedge\overline{\nabla}\!_A\,\epsilon\right).
\endalign
$$
Here $\overline{\nabla}\!_A\,\epsilon
=\bar\partial\epsilon+[A,\epsilon]$ is an infinitesimal
gauge transformation of a $\bar\partial$-connection.\qed
\enddemo

In the same way, the commutation relations of these Hamiltonians get
centrally extended:
$$
\{H_{\epsilon_1},H_{\epsilon_2}\}= H_{[\epsilon_1, \epsilon_2]}
+c(\epsilon_1, \epsilon_2)~,
$$
by the following cocycle (cf., eq.({2})):

\proclaim{Proposition 3.4}
$$
c(\epsilon_1, \epsilon_2)=2\pi i \int_P \operatorname{res}_P\,\sigma
\wedge
\operatorname{tr}(\epsilon_1\bar\partial\epsilon_2)~.
\tag{5}
$$
\endproclaim

Let $\hat{\goth g}^{(S,\sigma)}$ denote the Lie algebra of gauge
transformations on $S$ centrally extended by the cocycle (5),
and $\hat{G}^{(S,\sigma)}$ be the corresponding group (cf., \cite{FK})
\footnote{~It is not important which of the possible central extensions
of the group $G^S$ will be taken here since it is only the Lie algebra
that matters. However, if one will be able to consider the quantum
counterpart of
the Poisson geometry a difference may appear.}.
The infinite-dimensional space $\hat{\goth g}^{(S,\sigma)}$ is the
space of pairs $(\epsilon, z)$, where $\epsilon$ is a $\goth g$-valued
function on the surface $S$ and $z$ is a complex number.

We will take the following space of triples  $(F,C,x)$ as the space
$(\hat{\goth g}^{(S,\sigma)})^*$
dual to $\hat{\goth g}^{(S,\sigma)}$.
 Here $F$ is
a $(\text{End}\,E)$-valued (0,2)-form on $S$, ~$C$ is a
$(\text{End}\,E)$-valued (0,1)-form on
the ``boundary" (i.e., polar set) $P$ of $S$, and
$x$ is a complex number.  The nondegenerate
pairing $<\,,>$ between the spaces $\hat{\goth g}^{(S,\sigma)}$ and
$(\hat{\goth g}^{(S,\sigma)})^*$ is the following:
$$
<(F,C,x), (\epsilon, z)>=
\int_S \sigma \wedge \operatorname{tr}(\epsilon F)- 2\pi i \int_P
\operatorname{res}_P\,\sigma \wedge \operatorname{tr}(\epsilon C) +zx~.
$$

Let us consider the action of $\hat{G}^{(S,\sigma)}$ on $\Cal A^S$
generated
by the action of $G^S$. That is to say, the center of
$\hat{G}^{(S,\sigma)}$
acts trivially.

\proclaim{Proposition 3.5}
The centrally extended group of gauge transformations, $\hat
G^{(S,\sigma)}$,
acts on $\Cal A^S$  in a Hamiltonian  way.
The momentum map for the action of the corresponding gauge algebra
$\hat{\goth g}^{(S,\sigma)}$ is the mapping
$\Cal A^S\to (\hat{\goth g}^{(S,\sigma)})^*$ given by the
$(0,2)$-curvature
and by the restriction of the $\bar\partial$-connection form to the
``boundary":
$$
A\mapsto (\bar\partial A+A\wedge A, ~~A|_P, ~~1)~.
$$
\endproclaim

Let us denote, as before, by $G^S_P$ the group of gauge
transformations on $S$ based on $P$:
$G^S_P=\{g\in G^S~|~g|_P =\operatorname{id}\}$, and by
${\goth g}^S_P$ the corresponding Lie algebra.

Modifying slightly the last proposition, one obtains the following

\proclaim{Corollary 3.6}
The group $G^S_P$ acts on $\Cal A^S$ in a Hamiltonian way.
The momentum map for the action of the corresponding Lie algebra
${\goth g}^S_P$ is the mapping
$\Cal A^S\to (\goth g^S_P)^*$ given by the curvature:
$$
A\mapsto \bar\partial A+A\wedge A~.
$$
\endproclaim

\remark{Remark} Note that the group $G^\Sigma_\Gamma$ is not centrally
extended, but still $G^S_P\subset \hat G^{S,\sigma}$.
\endremark

Consider now the (holomorphic) Hamiltonian reduction
$\Cal A^S//G^S_P$ of the space of $\bar\partial$-connections
$\Cal A^S$ with respect to the group $G^S_P$. The result will be
the space of integrable
$\bar\partial$-connections in the bundle $E$ on $S$ modulo gauge
transformations from $G^S_P$. Such connections, with vanishing
$(0,2)$-form of the curvature tensor, are in one-to-one
correspondence
with {\it holomorphic structures} in the complex bundle $E$. Thus the
holomorphic Hamiltonian reduction leads us to the consideration of the
{\it space of all holomorphic structures in the bundle $E$ modulo gauge
equivalence  trivial on $P$}.
The corresponding quotient space
$\Cal M_{S,P}$,
which we consider only locally, near some of its smooth points, is,
by construction, an (infinite-dimensional) symplectic manifold.

On the other hand, let us consider now the space of
holomorphic structures in a smooth bundle on
a complex one-dimensional manifold by taking $P$ as such a manifold and
$E|_P$ as the bundle on $P$:
$$
\Cal C := \{\bar\partial_{\,0} +C\,|\, C\in
\Omega^{(0,1)}(P,\text{End}\,E|_P)\}~.
$$
Here $\bar\partial_{\,0}$ is also understood as the restriction to $P$
of our reference
holomorphic structure. The space $\Cal C$ of holomorphic
structures in a bundle on an elliptic curve (or a sum of such spaces if
$P$ consists of several disjoint components) is in fact an affine
subspace in a vector space dual to the Lie algebra $\hat\goth g^{P,\beta}$.
The latter is defined as the central extension of $\goth
g^P=\Gamma(P,\text{End}\,E|_P)$ by the cocycle
$$
c_\beta(\epsilon_1, \epsilon_2)=2\pi i \int_P \beta\,
\wedge
\operatorname{tr}(\epsilon_1\bar\partial\epsilon_2)~.
$$
This, of course, should be compared with eq.(5). In what follows we set
$\beta=\text{res}_P\,\sigma$ in which case
$\hat\goth g^{P,\beta}=\hat\goth g^{S,\sigma}/\goth g^S_P$.
The Lie
algebra $\hat\goth g^{P,\beta}$ (the two-dimensional current algebra,
in terminology of refs. \cite{EF, FK}, or ``double loop algebra'')
plays the role of the loop algebra, which appeared
in \S 1, while $\Cal C$ plays the role of the space of connections on a
circle. This point of view was suggested in refs. \cite{EF, FK}, where
more details can be found. We mention only the pairing which defines
$\Cal C$
as an affine subspace of the dual space to $\hat\goth g^{P,\beta}$:
$$
\align
<&(\epsilon,z)\,,(C,x)>=  2\pi i \int_P
\beta \wedge \operatorname{tr}(\epsilon C) +zx\,,~~~~~~~~\text{where}\\
&(\epsilon,z)\in\hat\goth g^{P,\beta}\,,(\bar\partial_{\,0}+C)\in\Cal
C,~~~\text{and}\,~~x,z\in\Bbb C ~.
\endalign
$$
Thus, the pairs $(C,1)$ define an affine subspace in
$(\hat\goth g^{P,\beta})^*$. We shall consider the standard
Kirillov--Kostant
Poisson structure on $(\hat\goth g^{P,\beta})^*$. Its symplectic leaves
are, as always, the coadjoint orbits, which in our case correspond to
isomorphism classes   of holomorphic bundles on $P$.

\proclaim{Proposition 3.7}
The mapping from the space
$\Cal M_{S,P}$ to the coadjoint representation space,
$(\hat\goth g^{P,\beta})^*$,
sending an integrable
$\bar\partial$-connection on the surface $S$ to its
restriction to $P$ is a Poisson mapping.
\endproclaim

\demo{Proof} This mapping is essentially the momentum mapping for the
action of gauge transformations on the boundary. \qed
\enddemo

Now we consider the quotient of the space $\Cal M_{S,P}$ by the whole
group $\hat G^{S,\sigma}$ of centrally extended gauge transformations.
The latter group acts on $\Cal M_{S,P}$ since gauge
transformations equal to the identity on $P$ form a
normal subgroup $G^S_P$ in $\hat G^{S,\sigma}$.

The quotient space
$$
\{\bar\partial +A\in \Cal A^S\,|\,\bar\partial A+A\wedge A=0\}/
\hat G^{(S,\sigma)}
$$
represents the set of isomorphism classes of holomorphic bundles on $S$
(corresponding to a given underlying topological bundle $E$). Then, by
construction, the local smooth moduli space $\Cal M_S$ of holomorphic
bundles on $S$ is a finite-dimensional Poisson manifold. The symplectic
leaves in it are described in terms of coadjoint orbits in
$\left(\hat{\goth g}^{P,\beta}\right)^*$
as follows (ref. \cite{FKR}).

\proclaim{Proposition 3.8} The local moduli space $\Cal M_S$ of
holomorphic bundles possesses a (holomorphic) Poisson structure.  The
symplectic leaves of this structure are parameterized by the moduli of
their
restrictions to the anticanonical divisor $P\subset S$ (that is a
symplectic
leaf is singled out by fixing the isomorphism class of the restriction
to
the elliptic curve $P$, or the isomorphism classes of restrictions to
each curve if $P$ consists of several such curves).
\endproclaim

\remark{Remarks} {\it i)} With minor modifications the above proposition
holds for $P$ consisting of several components intersecting
transversally.
In the latter case the corresponding gauge group $\hat G^{P,\beta}$
is the current group on a punctured Riemann surface, described in
ref. \cite{FK}.

 {\it ii)} The above proposition should not be understood
as that the isomorphism
classes of bundles on $P$ can be taken arbitrary; rather they have to
satisfy
the condition that they arise as restrictions of bundles defined over
$S$.
\endremark

\bigskip

{\bf Another description of symplectic structure.}
Here we give an alternative description
of the symplectic structure on the symplectic
leaves mentioned in Proposition 3.8.

Let $F$ be a holomorphic
bundle on $S$ corresponding to a smooth point in $\Cal M_S$, e.g.,
$H^2(S,\text{End}\,F)=0$
(we assume here that the structure group is now $G=GL(n,\Bbb C)$).
Let $K_S$ denote the canonical line bundle of $S$,
so that by assumption $K_S^{\;-1}$ possesses a holomorphic section, say,
$\eta$ and $P$ is the divisor of zeros of $\eta$.  Denote the sheaves
of holomorphic sections of the bundles under consideration
by the same symbols as the bundles themselves.
We can write down the following exact sequence of sheaves:
$$
0\to\text{End}\,F\otimes K_S\to\text{End}\,F\to\text{End}\,F\otimes\Cal
O_P \to 0~,
$$
where the second map is given by the multiplication by
$\eta$, while $\Cal O_P$ is the structure sheaf of the submanifold $P$.
The
associated long exact sequence reads as follows:
$$
\align
0\to & H^0(S,\text{End}\,F)\to H^0(P,\text{End}\,F|_P)\to \\
& H^1(S,\text{End}\,F\otimes K_S)\to H^1(S,\text{End}\,F)\to
  H^1(P,\text{End}\,F|_P)\to  \\
& H^2(S,\text{End}\,F\otimes K_S)\to 0~, \tag{6}
\endalign
$$
where we have used Serre duality to set
$H^0(S,\text{End}\,F\otimes K_S)=(H^2(S,\text{End}\,F))^*=0$.\,
\footnote{~The cohomology groups $H^0(S,\text{End}\,F)$ and
$H^2(S,\text{End}\,F\otimes K_S)$ appearing also in eq.({6}) are
of
equal dimension by Serre duality and are isomorphic to $\Bbb C$ if we
further
require $F$ to be a {\it simple} bundle. The groups
$H^1(S,\text{End}\,F)$ and
$H^1(S,\text{End}\,F\otimes K_S)$ are also related by Serre duality.}

   The cohomology group $H^1(S,\text{End}\,F)$ represents the space of
infinitesimal deformations of $F$. It describes the tangent
space to
$\Cal M_S$ at the point $F$
\;\footnote{~A description of the Poisson bivector on $\Cal M_S$ in
these terms was given in refs. \cite{Mu, Bon, Bot}.},
while the space tangent to the
symplectic leaf $\Cal S_F$ at $F$ should correspond (according to
Proposition 3.8) to such deformations of $F$
which leave $F|_P$ unchanged.  Thus, that tangent space corresponds to
the
kernel of the restriction map $$ T_F\,\Cal S_F =
\text{ker}\left(\,H^1(S,\text{End}\,F)\to
H^1(P,\text{End}\,F|_P)\,\right),
$$
or, by the exactness of ({6}), to the quotient
$$
T_F\,\Cal S_F = H^1(S,\text{End}\,F\otimes K_S) / I~,
$$
where
$$
I = \text{im}\left(\,H^0(P,\text{End}\,F|_P)\to
H^1(S,\text{End}\,F\otimes K_S)\,\right)~.
$$
Now consider the pairing
$$
H^1(S,\text{End}\,F) \otimes H^1(S,\text{End}\,F\otimes K_S) \to
H^2(S,K_S) \cong \Bbb C ~,
$$
induced by the multiplication in cohomology and taking trace in
$\text{End}\,F$. This map can be restricted in the first
factor
to the subspace $T_F\,\Cal S_F \subset H^1(S,\text{End}\,F)$, which
gives us
$$
T_F\,\Cal S_F \otimes H^1(S,\text{End}\,F\otimes K_S) \to \Bbb C~.
$$
By considering the second factor, one observes that this pairing
vanishes on
the
subspace $T_F\,\Cal S_F \otimes I$, thus, descending to a map of
$T_F\,\Cal S_F \otimes \left(H^1(S,\text{End}\,F\otimes K_S)/I \right)$
or, finally,
$$
T_F\,\Cal S_F \otimes T_F\,\Cal S_F \to \Bbb C~.\tag{7}
$$
One can check now that this pairing is skew-symmetric and defines
a 2-form on $\Cal S_F$ which coincides with the symplectic structure
discussed in Proposition 3.8.

It is more difficult  to prove the closedness of the 2-form (7)
in the latter description than by means  of  the Hamiltonian reduction
discussed above.
Note also that the alternative construction also has an  analogue
for the moduli space of flat connections, which is recovered by substituting
locally constant sheaves instead of sheaves of holomorphic sections.

\medskip


{\eightpoint \remark{Acknowledgements}
We are grateful to A.\,Bondal, V.\,Fock, and  A.\,Levin for fruitful
discussions. The work of B.K. was partially supported by an Alfred P. Sloan
Research Fellowship, by the NSF grants DMS-9627782 and DMS-9304580 (IAS),
and by the NSERC grant OGP-0194132.
The work of A.R. was partially supported by the grants
RFBR-98-01-00344, INTAS-96-0482 and the grant \# 96-15-96455 for the
support of scientific schools.
\endremark}

\bigskip


\end
\Refs

\ref \key Ar \by V.I.\,Arnold
\paper  Arrangement of ovals of real plane algebraic curves,
 involutions
of smooth four-dimensional  manifolds, and on arithmetics of
integral-valued quadratic forms
\jour Func. Anal. and Appl.  \vol 5:3   \yr 1971 \pages  169--176
\moreref
\paper Remarks on eigenvalues and eigenvectors of Hermitian
matrices, Berry
phase, adiabatic connections and quantum Hall effect
\jour Selecta Mathematica, New Series \vol 1:1 \yr 1995 \pages 1--19
\endref


\ref \key AB
\by M.\,Atiyah and R.\,Bott
\paper  The Yang-Mills equations over a Riemann surface
\jour Philos. Trans. Roy. Soc. London  \vol 308 \yr 1982   \pages
523--615
\endref

\ref \key Bon \by A.\,Bondal
\paper unpublished
\endref

\ref \key Bot \by F.\,Bottacin
\paper Poisson structures on moduli spaces of sheaves over Poisson surfaces
\yr 1995  \vol 121:2\pages 421-436
\jour
\endref

\ref \key Do \by S.K.\,Donaldson  \pages 453--460
\paper Instantons and geometric invariant theory
\yr 1984 \vol 93:4
\jour Commun. Math. Phys.
\endref

\ref \key DT \by S.K.\,Donaldson and R.P.\,Thomas \pages 31--47
\paper Gauge theory in higher dimensions
\yr 1998 \vol
\jour The geometric universe (Oxford, 1996),
Oxford Univ. Press, Oxford
\endref

\ref \key EF \by P.I.\,Etingof and I.B.\,Frenkel \pages 429-444
\paper Central extensions of current groups in two dimensions
\yr 1994 \vol 165
\jour Commun. Math. Phys
\endref

\ref \key FK \by I.B.\,Frenkel and B.A.\,Khesin \pages 541-562
\paper Four-dimensional realization of two-dimensional
 current groups
\yr 1996 \vol 178
\jour Commun. Math. Phys
\endref

\ref \key FKT \by I.B.\,Frenkel, B.A.\,Khesin, and A.N.\,Todorov
\paper  Complex counterpart of CSW theory
and holomorphic linking
\yr   \vol
\jour in preparation
\endref


\ref \key FKR \by V.V.\,Fock, B.A.\,Khesin, and A.A.\,Rosly
\paper unpublished
\yr   \vol
\jour
\endref

\ref \key FR \by V.V.\,Fock and A.A.\,Rosly
\paper Poisson structures on moduli of flat connections on
Riemann surfaces and $r$-matrices
\jour preprint   \yr 1992
\moreref
\paper  Flat connections and polyubles
\yr 1993  \vol 95:2
\jour Theor. Math. Phys. \pages 526--534
\endref

\ref \key GS
\by P.\,Griffiths and W.\,Schmid
\paper Recent developments in Hodge theory: a discussion of techniques
and results
\jour in Discrete subgroups of Lie groups and applications to moduli,
Bombay Colloquium 1973, Oxford Univ. Press \yr 1975  \vol   \pages
31--127
\endref

\ref \key Kh
\by B.A.\,Khesin
\paper Informal complexification and Poisson structures on moduli spaces
\yr 1997  \vol 180
\jour AMS Transl., Ser. 2 \pages 147--155
\endref

\ref \key KR
\by B.A.\,Khesin and A.A.\,Rosly
\paper in preparation
\yr   \vol
\jour
\endref


\ref \key Ko
\by S.\,Kobayashi
\book Differential geometry of complex vector bundles
\publ Publ. Math. Soc. Japan (Iwanami Shoten and Princeton Univ. Press)
\yr 1987
\endref


\ref \key KN
\by P.B.~Kronheimer and H.~Nakajima \pages 263--307
\paper  Yang-Mills instantons on ALE gravitational instantons
\yr 1990 \vol 288:2
\jour Math. Ann.
\endref

\ref \key Le
\by J.\,Leray \pages 81-180
\paper La calcul diff\'erential et int\'egral sur une vari\'et\'e
analytique complexe (Probl\'em de Cauchy, III)
\yr 1959 \vol 87
\jour Bull. Soc. Math. France
\endref

\ref \key LMNS
\by  A.~Losev, G.~Moore, N.~Nekrasov, S.~Shatashvili\pages 130-145
\paper Four-dimensional avatars of two-dimensional RCFT
\yr 1996 \vol 46
\jour Nucl. Phys. Proc. Suppl.
\endref


\ref \key Mu
\by S.\,Mukai
\paper Symplectic structure of the moduli space of stable sheaves on
    an abelian or $K3$ surface
\jour Inv. Math. \vol 77 \yr 1984\pages 101--116
\endref


\ref \key Th \by  R.P.\,Thomas \pages 1-104
\paper Gauge theory on Calabi--Yau manifolds
\yr 1997 \vol
\jour Ph.D. thesis, Oxford
\moreref\pages 1-49
\paper A holomorphic Casson invariant
for Calabi-Yau 3-fold, and bundles on $K3$ fibrations
\yr 1998 \vol
\jour preprint IAS
\endref

\ref \key Tu
\by A.\,Tyurin
\paper Symplectic structures on the moduli spaces of vector bundles
on algebraic surfaces with $p_g>0$
\jour Math. Izvestia\vol 33 \yr 1987\pages 139-177
\endref

\ref \key We \by A.~Weinstein
\paper The local structure of Poisson manifolds
\yr 1983 \vol 18:3
\jour J. Diff. Geom.
\pages 523--557
\endref

\ref \key Wi \by E.\,Witten
\paper Chern--Simons gauge theory as a string theory.
\yr 1995
\jour The Floer memorial volume, Progr. Math. 133, \publ Birkh\"auser,
Basel
\pages 637-678
\endref

\vskip .5cm

B.K.: Department of Mathematics, University of Toronto,
Toronto, ON M5S 3G3, Canada
and \break
School of Mathematics, Institute for Advanced Study, Princeton,
NJ 08540, USA; \break
e-mail:\,\,khesin\@math.toronto.edu
\newline

A.R.: Institute of Theoretical and Experimental Physics,
B.Cheremushkinskaya 25, Moscow 117218, Russia;
e-mail:\,\,rosly\@vxitep.itep.ru \endRefs
